\documentclass{article}
\usepackage{maa-monthly}
\usepackage{color}

\renewcommand{\a}{\alpha}
\renewcommand{\b}{\beta}

\newcommand{\bea}{\begin{eqnarray}}
\newcommand{\eea}{\end{eqnarray}}
\newcommand{\f}[2]{\frac{#1}{#2}}
\newcommand{\eq}{&=&}
\newcommand{\nn}{\nonumber \\ }

\newcommand{\siki}[1]{Eq. (\ref{#1})}
\newcommand{\zu}[1]{Fig. \ref{#1}}

\newcommand{\fn}[1]{\footnote{#1}}

\theoremstyle{theorem}

\theoremstyle{definition}

\begin{document}

\markright{Submission}

\title{
Trigonometric Ratios Can Prove\\ the Pythagorean Theorem 
}
\author{Shoya Kise, 
Takesa Uehara, 
and Takashi Shinzato
}

\maketitle
\begin{abstract}
Recent interest in noncircular trigonometric proofs has underscored the need for alternative methodologies. Jackson and Johnson's 2024 study addresses a longstanding gap in the foundations of trigonometric proofs.
Inspired by the work of Jackson and Johnson \cite{Jackson2024}, we present three noncircular proofs of the Pythagorean theorem based on  trigonometric identities.
First, we establish the Pythagorean theorem via an isosceles triangle construction 
and the tangent double-angle formula. Second, we present an alternative proof utilizing an isosceles-triangle and the angle-bisector theorem. 
Third, we derive a novel trigonometric relation from the angle‑bisector theorem, thereby unifying and extending the two preceding approaches.
These approaches collectively demonstrate that the principal contribution of Jackson and Johnson lies in their strategic use of the double‑angle formula.
These proofs clarify the role of trigonometric identities independent of infinite series.
\end{abstract}

\noindent

\section{Pythagorean Theorem}
The Pythagorean theorem occupies a central position in geometry and analysis. Consider triangle $ABC$ right-angled at $C$, with sides $AB=c,BC=a$ and $CA=b$. Then
\bea
c^2\eq a^2+b^2,
\label{eq1}
\eea
is well known as the Pythagorean theorem \cite{Jackson2024,Shirali2023,Villiers2023,Luvcic2022,Loomis1968,Shene2023, Lengvarszky2013,Basu2016,Bhatnagar2024,Zimba2009}. Dividing both sides of \siki{eq1} by $c^2$ yields $1=\f{a^2}{c^2}+\f{b^2}{c^2}$.
Setting $\sin\theta=\f{b}{c}$ and $\cos\theta=\f{a}{c}$ 
for $\angle ABC=\theta$, it follows that \siki{eq1} is equivalent to 
the Pythagorean identity $1=\cos^2\theta+\sin^2\theta$.

Historically, because the trigonometric ratios  $\sin\theta$ and $\cos\theta$  are defined by side‑length ratios in a right triangle, any proof of the Pythagorean theorem based on these ratios has been deemed circular.
This perspective motivates our search for genuinely noncircular approach. 
For this problem, Jackson and Johnson \cite{Jackson2024} 
bypass this obstacle by combining a geometric construction with an infinite geometric series and the double-angle identity, marking a breakthrough in the field \cite{Jackson2024,Shirali2023,Villiers2023,Luvcic2022}. 
Jackson and Johnson innovatively combined geometry and analysis in their approach. 
Their key insight was that combining an infinite geometric series with a double-angle identity circumvents circular reasoning.
Nonetheless, earlier works have provided proofs based exclusively on infinite geometric series
\cite{Shene2023,Lengvarszky2013,Bhatnagar2024}. 
This raises the question: is the novelty of Jackson and Johnson's approach attributable primarily to trigonometric identities or to their use of infinite series?

To gain deeper insight into their method, which resolved a problem unsolved for nearly two millennia, 
we analyze whether their success derives from the infinite series, the double‑angle formula, or a combination thereof.
Although prior research has established proofs of the Pythagorean theorem via infinite series  \cite{Loomis1968,Shene2023,Lengvarszky2013,Bhatnagar2024}, relatively, few have focused exclusively on the double-angle formula to prove the theorem.
Specifically, by isolating the double‑angle identity, we aim to demonstrate that infinite series are not essential for non‑circular trigonometric proofs
of the Pythagorean theorem.

Hence, we investigate whether the Pythagorean theorem can be derived solely from the double‑angle formula, without invoking infinite series.
Furthermore, we attempt to prove the Pythagorean theorem by employing the angle‑bisector theorem in conjunction with the double‑angle approach.
In addition, we discuss the connection between these two proofs, derive a new relational formula for trigonometric ratios, and use it to prove the Pythagorean theorem.
We also compare our approaches with classical Euclidean proofs to highlight their novelty.

This paper is organized as follows: Section \ref{Sec2} presents the first proof via an isosceles triangle construction and the tangent double-angle formula.  Section \ref{Sec3} offers an alternative proof using the angle-bisector theorem.  In section \ref{Sec4}, we derive the double-angle formula itself from the angle-bisector theorem, and section \ref{Sec5} develops a new relation to obtain a third proof.  In the final section, we conclude our work and devote future works.

\begin{figure}[b]
\begin{center}
\includegraphics[width=0.9\hsize,angle=0]{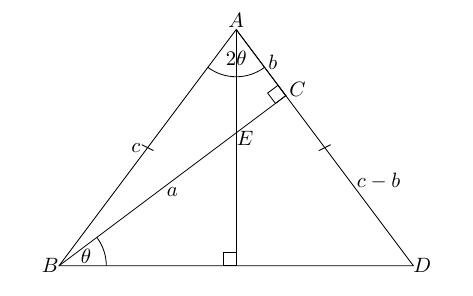}
\caption{\label{Fig1}Isosceles-Triangle Construction.}
\end{center}
\end{figure}

\section{First Proof}
\label{Sec2}
We begin by constructing an isosceles triangle that avoids reliance on series expansions.  Figure \ref{Fig1} thus 
illustrates the key angles and lengths in the isosceles construction, and provides intuitive geometric insight into the double-angle relationship. Let $ABD$ be an isosceles triangle 
with $AB=AD$ and $\angle BAD<\f{\pi}{2}$. Dropping the perpendicular from 
$B$ to $AD$ at $C$ yields a right triangle 
$ABC$ with $AB=c, BC=a$, and $CA=b$, hence $CD=AD-AC=c-b$. Setting $\angle BAD=2\theta$ gives
$\angle ADB=\angle ABD=
\f{\pi}{2}-\theta$. In the right triangle $DBC$, $\angle DBC=\f{\pi}{2}-\angle BDC=
\f{\pi}{2}-\angle ADB=
\theta$ is obtained.

Therefore, in right triangle $ABC$ (with apex angle $\angle BAD=2\theta$), 
$\tan2\theta=\f{BC}{CA}=\f{a}{b}$, 
and in right triangle $DBC$ (with apex angle $\angle DBC=\theta$), 
$\tan\theta=\f{CD}{BC}=\f{c-b}{a}$.
Substituting into the tangent double-angle formula
$\tan2\theta=\f{2\tan\theta}{1-\tan^2\theta}$
yields
\bea
\label{eq2}
\f{a}{b}\eq\f{2\left(\f{c-b}{a}\right)
}{1-\left(\f{c-b}{a}\right)^2},
\eea
then, rearranging \siki{eq2} by cross-multiplication immediately yields $c^2=a^2+b^2$. Thus, this construction provides an intuitive geometric interpretation and a noncircular proof of the Pythagorean theorem that relies solely on the tangent double-angle formula, without invoking infinite series.
In particular, it demonstrates that, between 
the two prior ideas-- (1) the infinite series and (2) the double-angle identity -- the latter is the essential ingredient.

Finally, two points should be noted here. First, since the above proof also covers the case 
$a=b$ (i.e.\ when triangle 
$ABC$ is an isosceles right triangle), it constitutes an extension of Jackson and Johnson's prior work. Next, 
it turns also out that this proof does not make use of the Pythagorean identity 
$1=\cos^2\theta+\sin^2\theta$.

\section{Second Proof}
\label{Sec3}
To validate the generality of our approach, we employ classical theorems in a novel manner.
In the previous section, we prove the Pythagorean theorem using the tangent double-angle formula. However, because we do not discuss whether that formula itself is derived from the Pythagorean identity $1=\cos^2\theta+\sin^2\theta$, 
we can not yet assert that this proof is entirely free of circularity.
Hence it is necessary to reexamine the proof by employing a different approach. Here we reuse the same setup in \zu{Fig1} but invoke the angle-bisector theorem to furnish a second proof.
Notably, our method reveals a deeper connection between angle-bisector properties and trigonometric relations.

Since we set 
$\angle BAD=2\theta$ in the first proof, let 
$E$  be the intersection of the bisector of 
$\angle BAD$ with side 
$BC$.
This similarity reveals a deeper structural link between angle bisectors and trigonometric ratios. Then right triangles 
$EAC$ and $DBC$ are similar, since 
$\angle EAC=\angle DBC=\theta$.
By the angle-bisector theorem,
$AC:AB=EC:EB$ and  $EC+EB=BC$, 
we compute
$EC=BC\f{AC}{AC+AB}=\f{ab}{b+c}$.
Hence, in the similar right triangles $EAC$ and $DBC$, $\tan\theta=\f{EC}{AC}=\f{\f{ab}{b+c}}{b}=\f{a}{b+c}$ and 
$\tan\theta=\f{DC}{BC}=\f{c-b}{a}$ are respectively obtained. Thus, equating the two expressions for 
$\tan\theta$ immediately gives
$c^2=a^2+b^2$.

In Section 4 we summarize how the tangent double-angle formula follows from the angle-bisector theorem and basic geometry, without ever invoking 
 the Pythagorean identity $1=\cos^2\theta+\sin^2\theta$. Moreover, a comparison of both proofs highlights the versatility of elementary trigonometric techniques. We now turn our attention to the implications of the unified identity.

\begin{figure}[b]
\begin{center}
\includegraphics[width=0.9\hsize,angle=0]{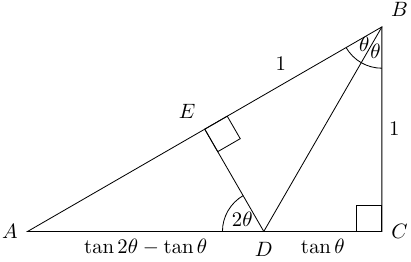}
\caption{\label{Fig2}Angle bisector theorem and double-angle formula.}
\end{center}
\end{figure}

\begin{figure}[b]
\begin{center}
\includegraphics[width=0.9\hsize,angle=0]{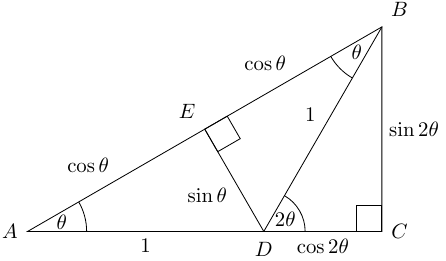}
\caption{\label{Fig3}Double angle formula.}
\end{center}
\end{figure}

\section{Derivation of double angle formula}
\label{Sec4}
Having established two independent proofs, we now derive standard identities from first principles.  In this section, we here derive the tangent double-angle formula 
$\tan2\theta=\f{2\tan\theta}{1-\tan^2\theta}$ from
the angle-bisector theorem,  without invoking the Pythagorean identity 
$1=\cos^2\theta+\sin^2\theta$. As shown in \zu{Fig2}, 
consider 
right triangle $ABC$ with $\angle ACB=\f{\pi}{2}$. 
We set $BC=1$ and let $D$ be the intersection of the bisector of $\angle CBA=2\theta$ with side $AC$.  From right triangle $BCD$, one finds  
$CD=\tan\theta$. 
We also set $AD=x$ and $AB=y$.  By the angle-bisector theorem $AB:BC=AD:CD$, hence  
$x=y\tan\theta$. 
Since $BC=1$, 
rewriting $\tan2\theta$ in terms of side lengths yields that 
\bea
\tan2\theta\eq\f{AC}{BC}=(y+1)\tan\theta,
\label{eq6}
\eea
where $AC=AD+CD=x+\tan\theta$ and $x=y\tan\theta$ are already used.
To complete the derivation, we use congruence arguments as follows. Next, drop the perpendicular from $D$ to side $AB$ and denotes its foot by $E$.  Since $\angle DBC=\angle DBE=\theta$, the right triangles $BCD$ and $BED$ are congruent by the hypotenuse angle condition; hence
$BC=BE=1$ and $CD=ED=\tan\theta$. 
Therefore, the right triangles $AED$ and $ACB$ are similar, implying $\angle ADE=2\theta$ and  
\bea
\label{eq7}
\tan2\theta\eq\f{AE}{ED}=\f{y-1}{\tan\theta},
\eea
where $AE=AB-BE=y-1$ is used.  Therefore, eliminating $y$ between  \siki{eq6} and \siki{eq7} yields the standard double-angle formula  
$\tan2\theta=\f{2\tan\theta}{1-\tan^2\theta}$.

Importantly, our derivation does not depend on the Pythagorean identity. 
Note that, although we use the definition of tangent, e.g.\ $\tan2\theta=\f{AC}{BC}=\f{AE}{ED}$, we do not rely on the Pythagorean identity $1=\cos^2\theta+\sin^2\theta$. Thus, one can derive the relation between $\tan2\theta$ and $\tan\theta$ without circularity (see Appendix \ref{app1}), so it is erroneous to claim that all trigonometric proofs of the Pythagorean theorem are circular by nature; a careful, case-by-case analysis yields deeper insight into each proof\fn{
We also derive the double-angle formulas for sine and cosine by constructing an appropriate right triangle and applying the law of sines and cosines. As shown in \zu{Fig3}, we prepare a right triangle $ABC$ with $\angle ACB=\f{\pi}{2}, \angle BAC=\theta\left(<\f{\pi}{4}\right)$. Next, place a point $D$ on side $BC$ so that $\angle ABD=\theta$, and for simplicity assume $AD=1$.  
Triangle $ABD$ is isosceles with $\angle ABD=\angle BAD=\theta$ and $AB=AD$. Since $\angle BDC=2\theta$, we obtain  
$BC=\sin2\theta$ and $CD=\cos2\theta$.
Let $E$ be the foot of the perpendicular from $D$ to side $AB$. Triangle $ADE$ is a right triangle with $\angle DEA=\f{\pi}{2}, \angle DAE=\theta$ and $ AD=1$. Hence  
$AE=\cos\theta$ and $AB=2AE=2\cos\theta$.
Therefore, the sine of angle $\angle BAC=\theta$ is  
$\sin\theta=\f{BC}{AB}=\f{\sin2\theta}{2\cos\theta}$ which yields the double-angle formula for sine, 
$\sin2\theta=2\sin\theta\cos\theta$.
Similarly, the cosine of $\angle BAC=\theta$ is  
$\cos\theta=\f{AC}{AB}=\f{1+\cos2\theta}{2\cos\theta}$ giving the double-angle formula for cosine, 
$\cos2\theta=2\cos^2\theta-1$, 
 where $AC=AD+CD=1+\cos2\theta$ is already used.
}.

\section{Third Proof}
\label{Sec5}
We extend our analysis by exploring alternative trigonometric relations. 
In this section, we employ the right triangle $ABC$ setup from Section \ref{Sec4} to derive 
additional trigonometric identities and present a novel proof of the Pythagorean theorem. One recalls that 
$BC=1$ and $CD=ED=\tan\theta$. 
Since 
$\angle ABC=2\theta$ and $BC=1$,  
we have $\tan2\theta=\f{AC}{BC}$, hence $AC=\tan2\theta$.  Similarly,  
$\cos2\theta=\f{BC}{AB}$ implies $AB=\f{1}{\cos2\theta}$.  
Next, the area $S$ of triangle $ABD$ can be expressed as  
$S=\f{1}{2}BC\cdot AD=\f{1}{2}AB\cdot ED$. 
Since 
$AD=AC-CD=\tan2\theta-\tan\theta$, one obtains 
$\f{\tan\theta}{\cos2\theta}=\tan2\theta-\tan\theta$ and 
derives  
$\tan\theta=\f{\sin2\theta}{1+\cos2\theta}$.

Moreover, we apply this identity to the configurations introduced in Sections \ref{Sec2} and \ref{Sec3}. 
In the right triangle $ABC$ 
with $\angle ACB=\f{\pi}{2}$, 
 $\angle BAC=2\theta$ and sides 
$AB=c,BC=a,CA=b$, we have $\sin2\theta=\f{BC}{AB}=\f{a}{c}$ and $\cos2\theta=\f{CA}{AB}=\f{b}{c}$. 
In the right triangle $BDC$ 
with 
$\angle BCD=\f{\pi}{2}$, $\angle DBC=\theta$, $BC=a$ and $CD=c-b$, 
one takes $\tan\theta=\f{CD}{BC}=\f{c-b}{a}$. 
Substituting into $\tan\theta=\f{\sin2\theta}{1+\cos2\theta}$ immediately yields $c^2=a^2+b^2$.  
Thus, this section demonstrates that, beyond
 the standard double-angle formula  
$\tan2\theta=\f{2\tan\theta}{1-\tan^2\theta}$ 
the Pythagorean theorem can also be proved using alternative trigonometric relations.

\section{Conclusion}
In this work, we have highlighted the primacy of double‑angle identities. 
In trigonometric proofs of the Pythagorean theorem, we have focused on two main approaches: (1) the infinite geometric series argument and (2) double-angle identities, following Jackson and Johnson's prior work. 
Our results demonstrate that the essential ingredient is the tangent double-angle identities, rather than infinite series constructions.
Specifically, we prove the Pythagorean theorem by employing the double-angle formula for tangent,  
$\tan2\theta=\f{2\tan\theta}{1-\tan^2\theta}$, and
we show that the angle-bisector theorem yields an equivalent derivation.
Moreover, we derive the  trigonometric relation  
$\tan\theta=\f{\sin2\theta}{1+\cos2\theta}$ and use it to prove the Pythagorean theorem.

Although this study and prior works have provided proofs based on the double-angle formulas, as future works, applying these double-angle formulas to novel geometric configurations may yield additional noncircular proofs of the Pythagorean theorem. Furthermore, to develop new trigonometric proofs of the Pythagorean theorem, it will be necessary to integrate triple-angle and half-angle formulas, as well as sum-to-product and product-to-sum identities, in a manner that avoids circular reasoning.
Further exploration of half‑angle formulas may yield additional insights.

\section*{Acknowledgments}
We thank T. Doi and G. Igusa for fruitful and detailed discussions.
\bibliographystyle{amsalpha}
\bibliography{sample20250418}

\providecommand{\bysame}{\leavevmode\hbox to3em{\hrulefill}\thinspace}
\providecommand{\MR}{\relax\ifhmode\unskip\space\fi MR }
\providecommand{\MRhref}[2]{%
  \href{http://www.ams.org/mathscinet-getitem?mr=#1}{#2}
}
\providecommand{\href}[2]{#2}
\begin{thebibliography}{Zim09}

\bibitem[Ba24]{Bhatnagar2024}
Gaurav Bhatnagar and Sagar~Shrivastava and, \emph{An uncountable number of
  proofs of the pythagorean theorem}, Mathematics Magazine \textbf{97} (2024),
  no.~5, 536--544.

\bibitem[Bas16]{Basu2016}
Kaushik Basu, \emph{A new and rather long proof of the pythagorean theorem by
  way of a proposition on isosceles triangles}, The College Mathematics Journal
  \textbf{47} (2016), no.~5, pp. 356--360.

\bibitem[dV23]{Villiers2023}
Michael de~Villiers, \emph{Is a trigonometric proof possible for the theorem of
  pythagoras?}, Learning and Teaching Mathematics \textbf{2023} (2023), no.~34,
  22--27.

\bibitem[JJ24]{Jackson2024}
Ne'Kiya Jackson and Calcea Johnson, \emph{Five or ten new proofs of the
  pythagorean theorem}, The American Mathematical Monthly \textbf{131} (2024),
  no.~9, 739--752.

\bibitem[Len13]{Lengvarszky2013}
Zsolt Lengv$\acute{\rm a}$rszky, \emph{Proving the pythagorean theorem via
  infinite dissections}, The American Mathematical Monthly \textbf{120} (2013),
  no.~8, 751--753.

\bibitem[Loo68]{Loomis1968}
Elisha~Scott Loomis, \emph{The pythagorean proposition : its demonstrations
  analyzed and classified, and bibliography of sources for data of the four
  kinds of "proofs"}, Classics in mathematics education, National Council of
  Teachers of Mathematics, 1968.

\bibitem[Lu{\v{c}}22]{Luvcic2022}
Zoran Lu{\v{c}}i{\'c}, \emph{Who proved pythagoras's theorem?}, The
  Mathematical Intelligencer \textbf{44} (2022), no.~4, 373--381.

\bibitem[She23]{Shene2023}
Ching-Kuang Shene, \emph{A new approach to proving the pythagorean theorem}.

\bibitem[Shi23]{Shirali2023}
Shailesh Shirali, \emph{Two new proofs of the pythagorean theorem-part i}, At
  Right Angles (2023), 7--12.

\bibitem[Zim09]{Zimba2009}
Jason Zimba, \emph{On the possibility of trigonometric proofs of the
  pythagorean theorem}, Forum geometricorum, vol.~9, 2009, pp.~275--278.

\end{thebibliography}

\begin{figure}[b]
\begin{center}
\includegraphics[width=0.9\hsize,angle=0]{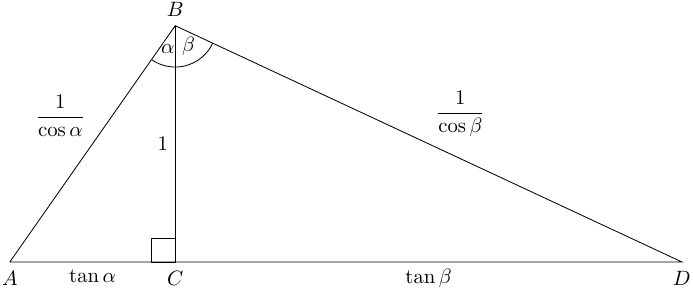}
\caption{Geometric proof of the addition angles formula.\label{Fig4}}
\end{center}
\end{figure}

\begin{figure}[t]
\begin{center}
\includegraphics[width=0.9\hsize,angle=0]{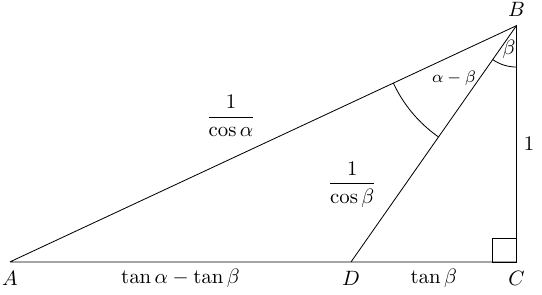}
\caption{Geometric proof of the subtraction angles formula.\label{Fig5}}
\end{center}
\end{figure}

\appendix

\section{Appendix}\label{app1}
For completeness, we derive the sine and cosine addition formulas geometrically. In the first part of this appendix, we discuss the addition angles formula for sine and cosine. 
As shown in \zu{Fig4}, first, let triangle $ABD$ satisfy $\angle BAD<\f{\pi}{2}$ and $\angle BDA<\f{\pi}{2}$, and let $C$ be the foot of the perpendicular from $B$ to side $AD$.  
This construction yields two right triangles: $ABC$ with $\angle ACB=\f{\pi}{2}$ and $DBC$ with $\angle DCB=\f{\pi}{2}$.  
To simplify the discussion, we set  
$BC=1,\angle ABC=\a$ and $\angle DBC=\b$, respectively. 
Then $AC=\tan\a, AB=\f{1}{\cos\a}, CD=\tan\b$ and $BD=\f{1}{\cos\b}$ are obtained.
Applying the law of sines to triangle $ABD$ yields
$\f{AD}{\sin(\a+\b)}=\f{AB}{\sin\left(\f{\pi}{2}-\b\right)}
=\f{BD}{\sin\left(\f{\pi}{2}-\a\right)}
$, where
$\angle ABD=\angle ABC+\angle DBC=\a+\b, AD=AC+CD=\tan\a+\tan\b$ and $\sin\left(\f{\pi}{2}-\a\right)=\cos\a$ are used. 
A straightforward rearrangement then gives the angle addition formula for sine $\sin(\a+\b)=\sin\a\cos\b+\cos\a\sin\b$.
Next, applying the law of cosines in triangle $ABD$ at the angle $\angle ABD=\a+\b$ yields
\bea
\cos(\a+\b)\eq\f{AB^2+BD^2-AD^2}{2AB\cdot BD}\nn
\eq\f{\cos\a\cos\b}{2}\left(\f{1}{\cos^2\a}+\f{1}{\cos^2\b}-(\tan\a+\tan\b)^2\right)\nn
\eq\cos\a\cos\b-\sin\a\sin\b,
\eea
where $\f{1}{\cos^2\a}=1+\tan^2\a$ is used\fn{
As shown in \zu{Fig6}, we prove the trigonometric identity $\f{1}{\cos^2\a}=1+\tan^2\a$ without invoking 
the Pythagorean identity $1=\cos^2\alpha+\sin^2\alpha$.  In right triangle $ABD$ with $\angle BAD=\f{\pi}{2}$, let $C$ be the foot of the perpendicular from $A$ to $BD$.  If $\angle ABC=\a$ and $BC=1$ are set, then $AC=\tan\a$ and $AB=\f{1}{\cos\a}$ are obtained.  Since $\angle DAC=\a$ and $\angle DCA=\f{\pi}{2}$, triangles $ABC$ and $DAC$ are similar, so
$CD=AC\tan\a=\tan^2\a$ and $BD=BC+CD=1+\tan^2\a$, are derived. Thus, we have
$\cos\a=\f{AB}{BD}=\f{\f{1}{\cos\a}}{1+\tan^2\a}$
whence $\f{1}{\cos^2\a}=1+\tan^2\a$ is derived.
}.
Thus, we also derive the angle addition formula for cosine.

\begin{figure}[b]
\begin{center}
\includegraphics[width=0.9\hsize,angle=0]{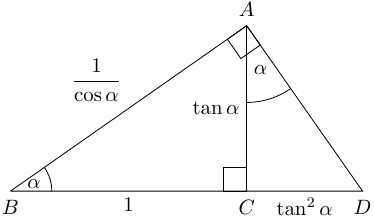}
\caption{$\f{1}{\cos^2\a}=1+\tan^2\a$.\label{Fig6}}
\end{center}
\end{figure}

\begin{figure}[t]
\begin{center}
\includegraphics[width=0.9\hsize,angle=0]{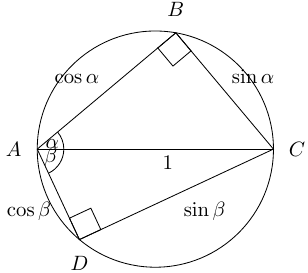}
\caption{Geometric derivation of the Pythagorean identity, $1=\cos^2\theta+\sin^2\theta$.\label{Fig7}}
\end{center}
\end{figure}

Furthermore, we derive the angle subtraction formulas. As shown in \zu{Fig6}, 
consider the right triangle $ABC$ with $\angle ACB=\f{\pi}{2}$ and $\angle ABC=\a$, and normalize $BC=1$.  Moreover, one puts $D$ on $AC$ so that $\angle DBC=\b(<\a)$.  Then
$AB=\f{1}{\cos\a}, BD=\f{1}{\cos\b}, AC=\tan\a$ and $CD=\tan\b$ are derived.
The law of sines in triangle $ABD$ yields
$\f{AD}{\sin(\a-\b)}=\f{BD}{\sin\left(\f{\pi}{2}-\a\right)}=\f{AB}{\sin\left(\f{\pi}{2}+\b\right)}$, by using $\angle ABD=\angle ABC-\angle DBC=\a-\b$ and $AD=AC-CD=\tan\a-\tan\b$, we find the angle subtraction formula for sine $\sin(\a-\b)=\sin\a\cos\b-\cos\a\sin\b$.

Finally, applying the law of cosines in triangle $ABD$ at $\angle ABD=\a-\b$ gives 
\bea
\cos(\a-\b)\eq\f{AB^2+BD^2-AD^2}{2AB\cdot BD}\nn
\eq\f{\cos\a\cos\b}{2}\left(
\f{1}{\cos^2\a}+\f{1}{\cos^2\b}
-(\tan\a-\tan\b)^2
\right)\nn
\eq\cos\a\cos\b+\sin\a\sin\b,
\eea
where $\f{1}{\cos^2\a}=1+\tan^2\a$ is employed.
Thus, the angle addition and subtraction formulas for sine and cosine can be derived without relying on the Pythagorean identity\fn{
As shown in \zu{Fig7}, we also prove the Pythagorean identity directly using a cyclic quadrilateral and the angle addition formulas \cite{Zimba2009}.  We prepare inscribe quadrilateral $ABCD$ in a circle of radius $\f{1}{2}$, and  $AC=1$ as a diameter, so that $\angle ABC=\angle ADC=\f{\pi}{2}$.  We set $\angle BAC=\a$ and $\angle DAC=\b$.  In right triangles $ABC$ and $ADC$, one has
$AB=\cos\a, BC=\sin\a, AD=\cos\b$ and $DC=\sin\b$. Then, the area $S$ of cyclic quadrilateral $ABCD$ can be written two ways:
$S=\f{1}{2}AB\cdot BC+\f{1}{2}AD\cdot DC$ and 
$S=\f{1}{2}AB\cdot AD\sin(\a+\b)+\f{1}{2}BC\cdot DC\sin(\pi-\a-\b)$, where 
$\angle BAD=\angle BAC+\angle DAC=\a+\b, \angle BCD=\pi-\angle BAD=\pi-\a-\b$, the angle addition formula for sine, and $\sin(\pi-\theta)=\sin\theta$, are already employed. Thus, one finds
\bea
0\eq\sin\a\cos\a(\cos^2\b+\sin^2\b-1)+\sin\b\cos\b(\cos^2\a+\sin^2\a-1).
\label{eq5}
\eea
Since \siki{eq5} holds for any $\a,\b$, it turns out that
$1=\cos^2\theta+\sin^2\theta$.
},  $1=\cos^2\theta+\sin^2\theta$.

\section{Exercise}
The following exercise consolidates the trigonometric techniques discussed above. As shown in \zu{Fig8}, we prepare right triangle $ABC$ with 
$\angle ABD=\angle BAD=\angle DBF=\theta\left(<\f{\pi}{6}\right), \angle BDF=2\theta, \angle BFC=3\theta, \angle ACB=\angle AED=\f{\pi}{2}, AD=BD=1, AB=2\cos\theta, BC=\sin2\theta$ and $CD=\cos2\theta$. Solve the following:
\begin{enumerate}
\item Find the lengths $BF$ and $DF$. 
\item Show $\sin3\theta+\sin\theta=2\cos\theta\sin2\theta$.
\item From the setting of \zu{Fig1}, since $\sin2\theta=\f{a}{c}, \cos2\theta=\f{b}{c}$ and $\tan\theta=\f{c-b}{a}$ are already obtained, prove the Pythagorean theorem.
\end{enumerate}

\begin{figure}[t]
\begin{center}
\includegraphics[width=0.9\hsize,angle=0]{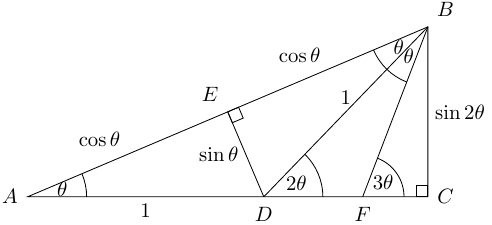}
\caption{Exercise.\label{Fig8}}
\end{center}
\end{figure}

\newpage
\section*{Shoya Kise}
 is a third-year junior high student at Kanagusuku Junior-High School in Okinawa, Japan. He often plays basketball with his best friends at a school sports club and enjoys watching mathematics lectures on YouTube. 

\section*{Takesa Uehara}
 is a first-year high school student at Kyuyo High School in Okinawa, Japan. He is interested in mathematics, physics, economics, and data analysis, plays a keyboard in his school music club; and enjoys online gaming with his best friends. 

\section*{Takashi Shinzato}
is a director of Okinawa Mathematical Academy, Okinawa, Japan, and  earned his PhD. in Science from  Institute of Science Tokyo, Tokyo, Japan in 2009. 
His research interests include statistical physics, machine learning, and financial engineering. E-mail: 
\href{mailto:takashi.shinzato@gmail.com}{{takashi.shinzato@gmail.com}}

\end{document}